\documentclass{elsarticle}

\usepackage[hidelinks]{hyperref}
\usepackage{amsmath}
\usepackage{amssymb}
\usepackage{amsthm}
\usepackage{tikz}

\newtheorem*{counterexample*}{The counterexample}

\newtheorem{innercustomthm}{Theorem}
\newenvironment{customthm}[1]
{\renewcommand\theinnercustomthm{#1}\innercustomthm}
{\endinnercustomthm}

\title{Note on clique polynomials and independent set polynomial of graphs}
\author{Hany Ibrahim}
\ead{hhasan@hs-mittweida.de}
\address{Applied university Mittweida, Depertment of Mathematics}

\begin{document}

\begin{abstract}
In 1994, Cornelis Hoede and Xueliang Li introduced the clique polynomial of a graph. Also, a theorem for the edge subgraph expansion for clique polynomials. In this note we present a counter example for it and explain which case it could be a valid theorem.
\end{abstract}

\maketitle

\section*{Introduction}
For notions, definitions and other reference, consult \cite{hoede1994clique}. We copy here the theorem under consideration as it is and we index the formulas in its proof.

\begin{customthm}{3.4} \cite{hoede1994clique}\label{TheTheorem}
	(Edge-subgraph expansion). Let $M$ be a nonempty subset of the edge set of a graph $G$. Then we have 
	\begin{equation*}
	C(G;x) = C(G-M;x) + \sum_{r=2}^{|M|}(-1)^{r} \sum_{\substack{S \subset M \\ |S|=(1/2)r(r-1) } } (r-1)x^{r}C^{*}(G_{S};x)
	\end{equation*}
	
	where
	\begin{equation*}
	C^{*}(G_{S};x) = 
		\begin{cases}
		C(G_{S};x)  & \text{if} \left< S \right> \text{is a clique of } G \text{,}\\
		0			& \text{otherwise.}
		\end{cases}
	\end{equation*}
\end{customthm}

\begin{proof}
	Let $H_{k}$ be the set of $k$-cliques of $G$. Then $|H_{k}| = a_{G}$. Let $M = {e_{1}, e_{2}, \cdots , e_{h}}$. Condition $Ci (i = 1,2, \dots , h)$, reads ‘contains edge $e_{i}$’. Then the number of $k$-cliques in $H_{k}$ that satisfy conditions $C_{i_{1}}, \cdots , C_{i_{t}}$, is given by
	\begin{align}
		N(C_{i_{1}},C_{i_{2}}, \cdots , C_{i_{t}}) =& a^{*}_{k - \left|V \left(\left<\left< e_{i_{1}}, e_{i_{2}},\cdots , e_{i_{t}}  \right>\right> \right)\right|} \left( G_{  \{e_{i_{1}}, e_{i_{2}}, \cdots , e_{i_{t}} \} } \right)  \nonumber \\
		=&   
		\begin{cases}
		a_{k - t^{*}} \left( G_{  \{e_{i_{1}}, e_{i_{2}}, \cdots , e_{i_{t}} \} } \right)  & \parbox[t]{.3\textwidth}{
if $\left<\left< e_{i_{1}}, e_{i_{2}},\cdots , e_{i_{t}}  \right>\right>$ is a $t^{*}$-clique of $G$,
	}
\\
			0	& \text{otherwise.}\\				
			\end{cases}
	\end{align}

	By the principle of inclusion and exclusion again we have the same formula as in the proof of Theorem $3.1$ ($M$ replaces $U$ and $e$’s replace $u$’s):
	
	\begin{align} \label{equation2}
	a_{k}(G-M) =& a_{k}(G) - \sum_{1 \leq i \leq h} a^{*}_{k - \left|V \left(\left<\left< e_{i}  \right>\right> \right)\right|} \left( G_{ \{ e_{i} \}} \right)  \nonumber \\
	+& \sum_{1 \leq i \leq j \leq h} a^{*}_{k - \left|V \left(\left<\left< e_{i} ,e_{j} \right>\right> \right)\right|} \left( G_{ \{ e_{i},e_{i} \}} \right)  \nonumber \\
	-& \sum_{1 \leq i \leq j \leq l \leq h} a^{*}_{k - \left|V \left(\left<\left< e_{i} ,e_{j},e_{l} \right>\right> \right)\right|} \left( G_{ \{ e_{i},e_{i},e_{l} \}} \right)  + \cdots . 
	\end{align}

The difficulty in this formula stems from the fact of the following sort. A $K_{6}$ can be induced by the vertices of $3$ independent edges, while $4,5,6, \cdots ,15$ edges may also induce a $K_{6}$. We want to rearrange the sums according to the orders of the graphs $\left<\left< e_{i_{1}},\cdots , e_{i_{t}}\right>\right>$. We observe: 
\begin{enumerate}
	\item If $\left< e_{i_{1}}, \cdots, e_{i_{t}}\right>$ is a clique, then $t= \binom{P}{2}$ for some $p$, i.e. $ |V(\left< e_{i_{1}}, \cdots , e_{i_{t}}\right>)| =p$. 
	\item A clique $K_{p}$, can have a spanning subgraph with $q$ edges for $ \lceil \frac{p}{2} \rceil \leq q \leq \binom{P}{2}$, where $\lceil x \rceil$ is the smallest integer not smaller than $x$. Contributions to $a_{p}(G - M)$ for a specific value $q$ stem from different terms in the expansion with a sign $(- l)^{q}$. 
	\item By Lemma 1.1 there are $f(p,q)$ spanning subgraphs with $q$ edges in $K_{p}$. These differ from the spanning subgraphs with $q$ edges in another graph $K_{p}$, as at least one vertex is different. Summing the contributions of edge sets that induce a specific $K_{p}$, gives 
	\begin{equation}
	\sum_{q=1}^{\binom{P}{2}} (-1)^{q}f(p,q)
	\end{equation}	
	which, by Lemma 1.2, is $(-1)^{p}p(p - 1)$. The numbers of $k$-cliques containing a specific $p$-clique is 
	
	{\small \begin{equation}
	a^{*}_{k - p} \left( G_{  \{e_{i_{1}}, e_{i_{2}}, \cdots , e_{i_{t}} \} } \right)  
	=
	\begin{cases}
	a_{k - p} \left( G_{  \{e_{i_{1}}, e_{i_{2}}, \cdots , e_{i_{\binom{P}{2}}} \} } \right)  & \parbox[t]{.28\textwidth}{
		if$\left< e_{i_{1}}, e_{i_{2}},\cdots , e_{i_{\binom{P}{2}}} \right>$ is a $p$-clique,}\\
	0	& \text{otherwise.}\\				
	\end{cases}
	\end{equation}}
\end{enumerate}

So, finally, we have 
	\begin{align}\label{equation3}
	a_{k}(G-M) =& a_{k}(G) - 1\sum_{1 \leq i \leq h} a^{*}_{k - 2} \left( G_{ \{ e_{i} \}} \right)  + 2 \sum_{1 \leq i < j \leq l \leq h} a^{*}_{k - 3} \left( G_{ \{ e_{i},e_{i},e_{l} \}} \right)  \nonumber \\
		&-3 \sum_{1 \leq i_{1} \leq i_{2} < \cdots < i_{6} \leq h } a^{*}_{k - 4} \left( G_{ \{ e_{i_{1}}, \cdots , e_{i_{6}} \}} \right)
	+ \cdots . 
	\end{align}
	
	By the same method as in Theorem 3.1, we obtain
	\begin{equation*}
	C(G;x) = C(G-M;x) + \sum_{r=2}^{|M|}(-1)^{r} \sum_{\substack{S \subset M \\ |S|=(1/2)r(r-1) } } (r-1)x^{r}C^{*}(G_{S};x)
	\end{equation*}
	Morever $C^{*}(G_{S};x)$ satisfies the condition in Theorem 3.4.
\end{proof}

\begin{counterexample*}
	The graph $G$ in figure \ref{counterexampleGraph} has $C(G;x)=1+3x+3x^{2}+x^{3}$. If we choose $M=\{ \{a,b\} , \{a,c\} \}$, then $C(G-M;x)=1+3x+x^{2}$ and the graph induced from the neighbors of $ \{a,b\}$ has  $C^{*}(G_{ \{a,b\} };x)=1+x$. By applying the formula in Theorem \ref{TheTheorem}, then we have
	\begin{align*}
		1+3x+3x^{2}+x^{3} \neq& (1+3x+x^{2}) + 2x^{2}(1+x) \\
						  \neq&  1+3x+3x^{2}+2x^{3}
	\end{align*}
\begin{figure}[ht!]
	\centering

		\begin{tikzpicture}
		\tikzset{vertex/.style = {shape=circle,draw}}
		
		\node[vertex] (a) at  (0,2) {a};
		\node[vertex] (b) at (-2,-2) {b};
		\node[vertex] (c) at (2,-2) {c};
		
		\draw (a) to (b) to (c) to (a);
		\end{tikzpicture}
		\caption{$G$}
		\label{counterexampleGraph}
\end{figure}

\begin{figure}[ht!]
	\centering
	
	\begin{tikzpicture}
	\tikzset{vertex/.style = {shape=circle,draw}}

	\node[vertex] (c) at (2,-2) {c};
	
	\end{tikzpicture}
	\caption{$G_{ \{a,b\} }$}
	\label{counterexampleGraph2}
\end{figure}

But, if we choose $M=\{ \{a,b\} , \{a,c\} , \{b,c\}\}$, then $C(G-M;x)=1+3x$ and the graph induced from the neighbors of $ \{a,b\}$ has  $C^{*}(G_{ \{a,b\} };x)=1+x$. $M$ induce and empty graph. By applying the formula in Theorem \ref{TheTheorem}, then we have
\begin{align*}
1+3x+3x^{2}+x^{3} =& (1+3x+x^{2}) + 3x^{2}(1+x) -2x^{3}\\
=&  1+3x+3x^{2}+x^{3}
\end{align*}
In general, the theorem is valid only when $M$ induces a clique. This is because, it was assumed that in equation \ref{equation2} there will be enough terms to be ordered and transformed to equation \ref{equation3}, but this is true only if $M$ induces a clique.
\end{counterexample*}


\end{document}